\documentclass[12pt]{amsart}

\usepackage{mathrsfs}
\usepackage{amsfonts}

\usepackage{amsmath}
\usepackage{amssymb}
\usepackage{multirow}
\usepackage{epsfig}

\usepackage{mathtools}
\usepackage{tensor}

\usepackage{a4wide}
\usepackage{pst-node}

\usepackage{tikz}
\usepackage{enumerate}

\newtheorem{propo}{Proposition}[section]

\newtheorem{lemma}[propo]{Lemma}

\newtheorem{theo}[propo]{Theorem}
\newtheorem{examp}[propo]{Example}

\newtheorem{rem}[propo]{Remark}

\newcommand{\bl}{\begin{lemma}}
\newcommand{\el}{\end{lemma}}

\def\d12{{_{12}}}

\def\Cay{{\rm Cay}}

\def\PG{{\rm PG}}

\def\qed{\ifmmode\square\else\nolinebreak\hfill
$\square$\fi\par\vskip12pt}

\renewcommand{\proof}{\par\noindent{\bf Proof.\ \ }}

\usepackage{indentfirst,latexsym,bm}

\def\GL{{\rm GL}}
\def\PGL{{\rm PGL}}
\def\AGL{{\rm AGL}}

\def\H{{\rm H}}
\def\K{{\rm K}}
\def\Sym{{\rm Sym}}

\def\Aut{{\rm Aut}}

\def\min{\mathrm{min}}
\def\diam{\mathrm{diam}}

\def\Aut{\mathrm{Aut}}

\def\Cay{{\rm Cay}}

\def\qed{\ifmmode\square\else\nolinebreak\hfill
$\square$\fi\par\vskip12pt}

\renewcommand{\proof}{\par\noindent{\bf Proof.\ \ }}

\begin{document}
\title{A classification of two-distance-transitive Cayley graphs over the generalized quaternion groups
}

\thanks{Supported by NSFC (12271524,12331013,12301453), NSF of Jiangxi (20224ACB201002) and NSF of Hunan (2022JJ30674)}

\author[W. Jin]{Wei Jin}
 \address{Wei Jin}
\address{School of Mathematics and Computational Science, Key Laboratory of Intelligent Computing and Information Processing of Ministry of Education\\
Xiangtan University\\
Xiangtan, Hunan, 411105, P.R.China}
\address{School of Statistics and Data Science\\
Jiangxi University of Finance and Economics\\
 Nanchang, Jiangxi, 330013, P.R.China}
\email{jinweipei82@163.com}
\author[P. S. Li]{Ping Shan Li}
\address{Ping Shan Li\\School of Mathematics and Statistics, Key Laboratory of Intelligent Computing and Information Processing of Ministry of Education\\
Xiangtan University\\
Xiangtan, Hunan, 411105, P.R.China}
\email{lips@xtu.edu.cn}
\author[L. Tan]{Li Tan  }
\address{Li Tan }
\address{School of Statistics and Data Science\\
Jiangxi University of Finance and Economics\\
 Nanchang, Jiangxi, 330013, P.R.China}
\email{tltanli@126.com}



\maketitle








\begin{abstract}

A  non-complete graph     is  \emph{$2$-distance-transitive} if, for $i=1,2$ and for any two vertex pairs $(u_1,v_1)$
and $(u_2,v_2)$ with the same distance $i$ in the graph, there exists an element of the graph automorphism group  that
maps  $(u_1,v_1)$ to $(u_2,v_2)$.
This    is  a generalization concept of the classical well-known  distance-transitive graphs.
In this paper, we completely determine the family of $2$-distance-transitive  Cayley graphs over the generalized quaternion groups.

\end{abstract}

\vspace{2mm}

\hspace{-17pt}{\bf Keywords:}  Cayley graph, 2-distance-transitive graph, generalized quaternion  group.

\hspace{-17pt}{\bf Math. Subj. Class.:} 05E18; 20B25

\section{ Introduction}
In this paper, all graphs are finite, simple, connected and undirected.
Let   $u$ and $v$ be two distinct vertices of a graph $\Gamma$.
Then the smallest  positive integer $n$ such that there is a path of length
$n$ from $u$ to $v$ is called the \emph{distance} from $u$ to $v$
and is denoted by $d_{\Gamma}(u, v)$.
Let $G$ be a subgroup of the automorphism group $\Aut(\Gamma)$ of $\Gamma$. Then  $\Gamma$ is said to be \emph{$(G,2)$-distance-transitive} if, for $i\in \{1,2\}$ and  for any two distinct vertex pairs $(u_1,v_1)$
and $(u_2,v_2)$ with $d_{\Gamma}(u_1,v_1)=d_{\Gamma}(u_2,v_2)=i$, there exists an element of  $G$  that maps $(u_1,v_1)$ to $(u_2,v_2)$.
If in particular $G=\Aut(\Gamma)$, then $\Gamma$ is simply said to be \emph{$2$-distance-transitive}.
The family  of $2$-distance-transitive graphs is a generalization  of the family of  classical well-known distance-transitive graphs.

The investigation of distance-transitive graphs is one in which the theoretical developments in algebraic graph theory can be traced back to the 1960s.
Efforts of many researchers (see \cite{BCN,DJLP-clique,HFZY-2022,JT-2016,Praeger-3,qdk-1}) have led almost to a complete classification of all finite distance-transitive graphs.
The systematic investigation of (locally) $2$-distance-transitive graphs was initiated recently.
Devillers, Giudici, Li and Praeger  ~\cite{DGLP-locdt-2012}
studied the class of locally $s$-distance-transitive graphs,
using the normal quotient strategy
developed for $s$-arc-transitive graphs in~\cite{Praeger-4}.
Corr, Schneider and the first author \cite{CJS-2dt} investigated the family  of $2$-distance-transitive  graphs of girth 4, and  determined  $2$-distance-transitive but not $2$-arc-transitive graphs of valency at most 5. After that,    the family of $2$-distance-transitive Cayley graphs over cyclic groups (circulants) is  identified    in \cite{CJL-2019}.
Recently,  a  classification   of 2-distance-transitive Cayley graphs over the dihedral groups was completed  in \cite{HFZY-2025,JT-2022}.



Let $T$ be  a finite group  and let $S$ be  a subset  of $T$ such that $1\notin
S$ and $S=S^{-1}$. Then the \emph{Cayley graph} $\Cay(T,S)$ of $T$ with
respect to $S$ is  the graph with vertex set $T$ and edge set
$\{\{g,sg\} \,|\,g\in T,s\in S\}$. In particular, the Cayley graph $\Cay(T,S)$ is
connected if and only if $T=\langle S\rangle$.
The group $R(T) = \{ \sigma_t|t\in T\}$ consists of right translations $\sigma_t : x \mapsto xt$
being a subgroup of the automorphism group $\Aut(\Gamma)$ and acting
regularly on the vertex set. We may identify $T$ with $R(T)$. C. Godsil
\cite[Lemma 2.1]{Godsil-1981} observed that
$N_{\Aut(\Gamma)}(T)=T:\Aut(T,S)$ where $\Aut(T,S)=\{\sigma \in
\Aut(T)|S^\sigma=S\}$. If $\Aut(\Gamma)=N_{\Aut(\Gamma)}(T)$, then the graph $\Gamma$ was called a \emph{normal Cayley graph} by Xu \cite{Xu-1} and such graphs have been studied  under
various additional conditions, see
\cite{DJLP-cayleyred1,DWX-1998,LX-2003,Pan-2014,Praeger-1999-cay}.

In this paper, we continue the study of the family of $2$-distance-transitive Cayley graphs, precisely, we are interested in
$2$-distance-transitive Cayley graphs over the generalized quaternion groups.

The  \emph{generalized quaternion group} or \emph{dicyclic group} of order $4n$  is   $Q_{4n}=\langle a,b\rangle$, where
$a^{2n}=1,b^2=a^{n}$ and $b^{-1}ab=a^{-1}$.
This group  is abelian if and only if $n=1$,  and the subgroup $\langle a\rangle\cong \mathbb{Z}_{2n}$ is a characteristic subgroup of $Q_{4n}$.

The following is a family of examples which are   $2$-distance-transitive Cayley graphs over the generalized quaternion groups.

\begin{examp}\label{gq-p-1}
{\rm Let $T= \langle a,b|a^{2n}=1,b^2=a^2,a^b=a^{-1}\rangle \cong Q_{4n}$ with $n\geq 2$, $S=T\setminus \langle a^n\rangle$ and $\Gamma:=\Cay(T,S)$.
Let $u=1_T$. Then $\Gamma_2(u)=\{a^n\}$, and $\{u\}\cup S\cup \Gamma_2(u)=T$.
Moreover,  $\Gamma$ has diameter 2, antipodal with each fold has 2 vertices, and so    $\Gamma \cong \K_{2n[2]}$.

Its automorphism group    $\Aut (\Gamma)=S_2\wr S_{2n}$ is transitive on both the set of
vertices  and  the set of arcs.
For each arc  $(u,v)$  of
$\Gamma$, we have  $|\Gamma_2(u)\cap \Gamma(v)|=1$,  and so $\Gamma$ is
$2$-distance-transitive. Since $\Gamma$ has girth 3 and non-complete, it follows that it  is not 2-arc-transitive.

}
\end{examp}

The graph in Figure \ref{fig1}  is $\K_{4[2]}$ which is a $2$-distance-transitive Cayley graph  over the generalized quaternion group $Q_8$.






\begin{figure}[t]\label{fig1}
\centering

\begin{tikzpicture}

\draw (7,4)-- (10,2); \draw (7,4)-- (10,3);

\draw (7,4)-- (7,0.7);

\draw (7,4)-- (8,0.7);

\draw (5,2)-- (7,4); \draw (5,3)-- (7,4);

\draw (8,4)-- (10,2); \draw (8,4)-- (10,3);

\draw (8,4)-- (7,0.7);

\draw (8,4)-- (8,0.7);

\draw (8,4)-- (5,2); \draw (5,3)-- (8,4);

\draw (5,3)-- (10,2); \draw (5,3)-- (10,3);

\draw (5,3)-- (7,0.7);

\draw (5,3)-- (8,0.7);

\draw (5,2)-- (10,2); \draw (5,2)-- (10,3);

\draw (5,2)-- (7,0.7);

\draw (5,2)-- (8,0.7);



\draw (7,0.8)-- (10,2); \draw (7,0.8)-- (10,3);

\draw (7,0.7)-- (8,0.7);

\draw (8,0.7)-- (10,2); \draw (8,0.7)-- (10,3);


\filldraw[black] (5,2) circle (2pt)  (5,3) circle (2pt);

\filldraw[black] (10,2) circle (2pt)  (10,3) circle (2pt);


\filldraw[black] (7,0.7) circle (2pt);


\filldraw[black] (8,0.7) circle (2pt);

\filldraw[black] (7,4) circle (2pt) (8,4) circle (2pt);


\end{tikzpicture}
\caption{$\K_{4[2]}$}
\end{figure}




Our main  theorem gives  a complete classification of the family of $2$-distance-transitive Cayley
graphs over the generalized quaternion  groups.



\begin{theo}\label{2dt-gq-theo1}
Let $\Gamma$ be a connected  $2$-distance-transitive Cayley graph over an order $4n$ generalized quaternion group, where $n\geq 2$.
Then $\Gamma$ is one of the following graphs:

 \begin{itemize}


\item[(1)]  $ \K_{x[y]}$ for some $x\geq 3,y\geq 2$ and $xy=4n$;

\item[(2)] $\K_{2n,2n}$;

\item[(3)]     $\K_{2n,2n}-2n\K_2$ where $n\geq 3$;

\item[(4)]   $B(PG(d,q))$ or $B'(PG(d,q))$, where $2n=(q^d-1)/(q-1)$, $d\geq 2$, and $q$ is a prime power;


\item[(5)]    $\K_{q+1}^{2d}$, where $q$ is an odd prime power,   $d\geq 2$   dividing $q-1$ and $2n=d(q+1)$;

\item[(6)]    $ X_1(4,q)$ where $q\equiv 3\pmod{4}$ and  $q=n-1$;

\item[(7)]  $\Gamma(d, q, r)$ where $r|q-1$;



\item[(8)] $X(2,2)$;

\item[(9)] $X'(3,2)$;




\item[(10)]   $ X_2(3)$.

\end{itemize}

\end{theo}

The definitions of the graphs arising in Theorem  \ref{2dt-gq-theo1} will be given in the next section.

This paper is organized  as follows.
After this introduction, we give, in Section 2,
some definitions on groups and graphs that we need and also prove some elementary  lemmas which will be used in the following analysis.
Let $\Gamma$ be a connected $(G,2)$-distance-transitive Cayley  graph over the generalized quaternion group    $T=\langle a,b|a^{2n}=1,a^n=b^2,a^b=a^{-1}\rangle$, where $n\geq 2$
and $T\leq G\leq \Aut(\Gamma)$.
 Suppose that $\Gamma\ncong \K_{m[f]}$ for any  $m\geq 3$ and $f\geq 2$.
 Let   $N$ be a  normal subgroup of $G$ maximal with respect to having at least $3$ orbits.
It is proved in Theorem \ref{bic-reduction-th2} that:
$\Gamma$ is a cyclic  cover of $\Gamma_N$  which
 is a $(G/N,s)$-distance-transitive Cayley graph over $T/N$, where $s=\min\{2,\diam(\Gamma_N)\}$, and $G/N$ is faithful and either quasiprimitive or bi-quasiprimitive on $V(\Gamma_N)$.
Moreover,   $N=\langle a^i\rangle$ where $i$ is a divisor of $2n$ and $i\neq 1,2n$,  and $T/N$ is known.
All the base graphs are determined, and  the cyclic covers of base graphs are found.


\section{Preliminaries}

In this section, we will give some definitions about groups and graphs that will be used in the
paper.

For a graph $\Gamma$, we use $V(\Gamma)$ and $\Aut(\Gamma)$ to denote
its \emph{vertex set}  and \emph{automorphism
group}, respectively. For the group theoretic terminology not defined here we refer the reader to \cite{Cameron-1,DM-1,Wielandt-book}.

\subsection{Groups and graphs}


A transitive permutation group $G\leqslant \Sym(\Omega)$ is said to be \emph{quasiprimitive}, if every non-trivial normal subgroup of $G$ is transitive on $\Omega$, while
$G$ is said to be \emph{bi-quasiprimitive} if every non-trivial normal subgroup of $G$ has at most two orbits on $\Omega$ and there exists one which has exactly two orbits on $\Omega$.
Quasiprimitivity  is a generalization
of primitivity as every normal subgroup of a primitive group is transitive,
but there exist quasiprimitive groups which are not primitive. For more information about quasiprimitive and bi-quasiprimitive permutation groups, refer to \cite{Praeger-1993-onanscott,Praeger-2,Praeger-2003-biq}.

A group extension $G = N.H$ of $N$ by $H$ is called   a \emph{central extension} if $N\leq Z(G)$, and
a central extension $G = N.H$ is referred to as a \emph{covering group} of $H$ if $G$ is perfect,
meaning the derived group $G'$ of $G$ equals $G$. Schur \cite{Schur-2009} showed  that a simple
group $N$ possesses a universal covering group $G$, that is, every covering group of $N$ is
a homomorphic image of $G$, and the center $Z(G)$ is called the \emph{Schur multiplier} of $N$,
denoted by $Mult(N)$.

We denote by  $\K_{m[b]}$  the \emph{complete multipartite graph} with $m$ parts, and each part has
$b$ vertices where $m\geq 3, b\geq 2$.

We call a graph with $n$ vertices  a {\it circulant} if it has an automorphism  that is an  $n$-cycle.
Thus a circulant is a Cayley graph over a cyclic group.

A {\it dihedral group} of order $2n$ is denoted by $D_{2n}$ and defined as $D_{2n}=\langle a,b|a^n=1,b^2=1,bab=a^{-1}\rangle$.
A  graph  is called a \emph{dihedrant} if it is a Cayley graph over a   dihedral group.

The    \emph{Hamming graph}  $\H(d,r)$  has
vertex set $\Delta^d=\{(x_1,x_2,\ldots,x_d)|x_i\in \Delta\}$, where
$\Delta=\{0,1,\ldots,r-1\}$, and  two vertices $v$ and $v'$ are adjacent if and only if they are different
in exactly one coordinate. The    Hamming graph  $\H(d,2)$ is  called a \emph{$d$-cube}, whose  valency is $d$.

For a graph $\Gamma$, its  \emph{diameter}   is the maximum distance occurring over all
pairs of vertices.
For  $u\in V(\Gamma)$ and each integer $i$ less than or equal to the diameter of $\Gamma$, we use $\Gamma_i(u)$ to denote the set of vertices at distance $i$ with vertex $u$ in $\Gamma$. Further
$\Gamma_1(u)$ is usually  denoted by $\Gamma(u)$.




A vertex triple $(u,v,w)$ of $\Gamma$  with $v$ adjacent to both $u$
and $w$ is called a \emph{$2$-arc} if $u\neq w$.
A $G$-arc-transitive
graph $\Gamma$   is said to be \emph{ $(G,2)$-arc-transitive}  if  $G$  is transitive on the set of
2-arcs.  Moreover, if $G=\Aut(\Gamma)$, then $G$ is usually omitted in the previous notation.
The first remarkable result about the class of finite $2$-arc-transitive graphs comes
from Tutte  \cite{Tutte-1,Tutte-2}.  Due to the above seminal work, this class  of graphs has
been studied extensively in the literature, see
\cite{ACMM-1996,IP-1,Li-abeliancay-2008,Praeger-1993-1,SDS-2019}.


Let   $d\geq 3$ be an integer and  $q$ be a prime power, and let $V$ be a $d$-dimensional vector space over the field $GF(q)$.
Set $B(\PG(d-1,q))$ to be the bipartite graph with vertices the 1-dimensional and $(d-1)$-dimensional subspaces of  $V$, and two subspaces are adjacent if and only if one is contained in the other. We denote $B'(\PG(d-1,q))$ to be the bipartite complement of $B(\PG(d-1,q))$, that is, the bipartite graph with the same vertex set as $B(\PG(d-1,q))$, but a 1-subspace and a $(d-1)$-subspace  are adjacent if and only if their intersection is the zero subspace.
Both of these two graphs have the same automorphism group $P\Gamma L(d,q):\mathbb{Z}_2$, refer to \cite[p.210]{CO-1987} and \cite{DMM-2008}.

The \emph{generalized Petersen graph} $GP(n,r)$  is the graph on the vertex set
$$\{u_0,u_1,\ldots,u_{n-1},v_0,v_1,\ldots,v_{n-1}\}$$
with the adjacencies:
$$u_i\sim u_{i+1}, \, v_i \sim v_{i+r}, \, u_i\sim v_i, \quad i=0,1,\ldots,n-1.$$
Hence each  generalized Petersen graph $GP(n,r)$ has valency 3.
By \cite{FGW-1971},
$GP(n,r)$ is arc-transitive if and only if   $(n,r)=(4,1)$, $(5, 2) $, $(8, 3)$, $ (10, 2),$ $ (10, 3), (12, 5) $  and $  (24, 5)$.

Let $\Gamma$ be a $G$-vertex-transitive graph. Let $\mathcal{B}=\{B_1,\ldots,B_m\}$ be a $G$-invariant partition of the vertex set, that is, for each $B_i$ and each $g\in G$, either $B_i^g\cap B_i= \varnothing $ or $B_i^g=B_i$.
Then the \emph{quotient graph}
$\Gamma_{\mathcal{B}}$ of $\Gamma$  is  the graph whose vertex set is the set $\mathcal{B}$,
such that two elements   $B_i,B_j \in \mathcal{B}$ are adjacent in $\Gamma_{\mathcal{B}}$ if
and only if there exist $x\in B_i$ and $ y\in B_j$ such that $x,y$ are
adjacent in $\Gamma$.
The graph    $\Gamma$ is  called a
\emph{cover} of $\Gamma_{\mathcal{B}}$ if,  for each edge $\{B_i,B_j\}$ of
$\Gamma_{\mathcal{B}}$ and $v\in B_i$, the vertex $v$ is adjacent to exactly one vertex in $B_j$; and further if $|B_i|=n$ and we want to emphasize this value, we  call
$\Gamma$  a
\emph{$n$-cover} of $\Gamma_{\mathcal{B}}$.
Moreover, if $\mathcal{B}$ is the set of orbits of an intransitive normal subgroup of $G$, then
$\Gamma$ is called a \emph{normal $r$-cover} of $\Gamma_{\mathcal{B}}$.

Gross and Tucker \cite{GT-1987} introduced a  combinatorial method  of  cover graphs through a voltage graph.
A \emph{graph homomorphism} $f$ from a graph $\widetilde{X}$ to a graph $X$ is a mapping from the vertex
set $V ( \widetilde{X} )$ to the vertex set $V (X)$ such that if $\{u, v\}\in  E( \widetilde{X} )$ then $\{f(u), f(v)\}\in  E(X)$.
When $f$ is surjective, $X$ is called a \emph{quotient of $\widetilde{X}$}. For $v\in  V (X)$, let $X(v)$ denote the
set of neighbours of $v$ in $X$. A \emph{covering projection} is defined as a graph homomorphism
$p : \widetilde{X} \mapsto X$ which is surjective and locally bijective, which means that the restriction
$p : \widetilde{X}(\widetilde{v}) \mapsto X(v)$ is a bijection, whenever $\widetilde{v}$ is a vertex of $\widetilde{X}$ such that $p(\widetilde{v}) = v \in V (X)$.
We call $X$ the \emph{base graph}, $\widetilde{X}$ a \emph{cover graph} of $X$, and the pre-images
$p^{-1}(v), v\in  V (X)$ the \emph{fibres}. A covering projection $p : \widetilde{X} \mapsto X$ is called \emph{regular} if there
exists a semi-regular subgroup $N$ of the group of automorphisms $\Aut( \widetilde{X} )$ of $\widetilde{X}$ such that
the quotient graph $\widetilde{X}_N$ (with vertices taken as the orbits of $N$ on $V( \widetilde{X} ))$ is isomorphic
to $X$. In that case we call $N$ the \emph{covering (transformation) group}. Particularly,
a regular covering projection $p : \widetilde{X} \mapsto X$  is said to be  a \emph{regular cyclic covering projection}  whenever the covering group $N$ is  cyclic, and $\widetilde{X}$ is said to be  a \emph{regular cyclic cover} of $X$.

The above properties can be exploited to construct regular cover graphs of a given
graph, as follows.

Let $X$ be a connected graph, and let $N$ be a finite group. Suppose $\psi : Arc(X) \mapsto N$ is
a function assigning a group element to each arc of $X$, such that $\psi(v, u) = (\psi(u, v))^{-1}$ for
every arc $(u, v)\in  Arc(X)$. Here $\psi$ is called a \emph{voltage assignment},  the values of $\psi$ are called \emph{voltages}, and $N$ is  the covering
group. In particular, $\psi$ is
called \emph{reduced} if the values of $\psi$ on a spanning tree are trivial (equal to the identity element
of $N$). We may construct a larger graph $Cov(X,\psi)=X\times_{\psi} N$, called the \emph{(derived) voltage graph}
(or \emph{cover graph}), with vertex set $V (X) \times N$ and adjacency defined by $(u, g) \sim  (v, h)$
if and only if $u \sim v$ and $h = g\psi (u, v)$. Moreover, a voltage assignment on arcs could be naturally extended to a voltage assignment on walks, such that if $W$ and $W'$ are two walks with the last vertex of $W$ being the same as the first vertex of $W'$, then in the voltage group $\psi(WW')=\psi(W)\psi(W')$.
For  a walk $W$ of $X$, if the voltage of $W$ is the identity of the voltage group, then we say that
 $W$ has a \emph{vanishing voltage}.

Let $q$ be a prime power such that $q\equiv 3\pmod{4}$ and let
  $S(q)$ be the set of all non-zero squares of $GF(q)$
and $N(q)$ be the set of  non-zero non-squares of $GF(q)$.
Let $Y=\K_{q+1}$, and   identify the vertices of $Y$ with the projective lines
$\PG(1, q)=GF(q)\cup \{\infty\}$ (see \cite[p.283,285]{DMW-1998}).
 We define $X_1(4, q)$ to be the 4-fold
cover $Cov(Y, f )$, where the voltage $\psi: Arc(Y)\mapsto \mathbb{Z}_4$ is defined with the following
rule:

\begin{center}
$\psi(x,y): = \left\{
\begin{array}{lll}
0, & \mbox{if $\infty \in\{x,y\} $};\\
1, & \mbox{if $ y-x \in S(q)$};\\
3, & \mbox{if $y-x \in N(q)$}.\end{array} \right.$
\end{center}

If $q=3$, then   $X_1(4, q)$ is the generalized Petersen graph $GP(8,3)$, and is also known as the M\"obius-Kantor graph.

We denote  the unique symmetric $(11,5,2)$-design by $H(11)$. Let  $A$ be the set of points  and $A'$ be the set of blocks of $H(11)$. Then the incidence graph $B(H(11))$ of this design has vertex set $A\cup A'$ and edge set
$\{xy|x\in A,y\in A',x\in y\}$. Let  $B'(H(11))$ be the incidence graph of the complementary design of $H(11)$. Then $B'(H(11))$ has vertex set $A\cup A'$ and edge set $\{xy|x\in A,y\in A',x\notin y\}$.

Let $q = r^l$ where $r$ is an odd prime  and let $GF(q)^*=\langle \theta \rangle$
be the multiplicative group of the field $GF(q)$ of order $q$. Let $\K_{q+1,q+1}-(q + 1)\K_2$ be
the graph that the  complete bipartite graph $\K_{q+1,q+1}$ minus a perfect matching whose vertex set is two copies of the
projective line $\PG(1, q)$, where the missing perfect matching consists of all pairs $[i,i'], i\in \PG(1, q)$. For any $d | q-1$ and $d\geq 2$, define a voltage graph
$\K_{q+1}^{2d}= (\K_{q+1,q+1}-(q + 1)\K_2)\times_f \mathbb{Z}_d$, where

$$f_{\infty',i}=f_{\infty,j'}=\overline{0}, i,j\neq \infty;$$
$$f_{i,j'}=\overline{h}, j-i=\theta^h, i,j\neq \infty.$$


Let $p$ be an odd prime and let $r$ be a positive   integer  dividing $p-1$. Let $A$ and $A'$ denote two disjoint copies of $\mathbb{Z}_p$ and denote the corresponding elements of $A$ and $A'$ by $i$ and $i'$, respectively.
Let $L(p,r)$ be the unique subgroup of the multiplicative group of $\mathbb{Z}_p$ of order $r$.   We define
two graphs,  $G(2p,r)$ and $G(2,p,r)$, with vertex set $A\cup A'$.  The graph $G(2p,r)$ has edge set $\{\{x,y'\}|x,y\in \mathbb{Z}_p,y-x\in L(p,r)\},$ while  the graph $G(2,p,r)$ (defined only for $r$ even) has  edge set $\{\{x,y\},\{x',y\},\{x,y'\},\{x',y'\}|x,y\in \mathbb{Z}_p,y-x\in L(p,r)\}.$
Note that $G(2,p,r)$ is a non-bipartite bicirculant of valency $2r$ as it contains a $p$-cycle.
Moreover, if $r=p-1$, then $G(2,p,r)$ is the graph $\K_{p[2]}$, and  is also the complement graph of a complete matching.

Let $\Gamma=\K_{5,5}-5\K_2$, where
$V(\Gamma)=\{1,2,3,4,5\}\cup\{1',2',3',4',5'\}$,
and $E(\Gamma)=\{ij'|i\neq j,i,j'\in V(\Gamma)\}$.

Define $X_2(3)=\Gamma \times_\psi \mathbb{Z}_3$, with the voltage assignment
$\psi:Arc(\Gamma)\mapsto \mathbb{Z}_3$ such that
$$\psi_{1,2'}=\psi_{1,3'}=\psi_{1,4'}=\psi_{1,5'}=\psi_{2,1'}=\psi_{2,3'}=\psi_{3,1'}$$
$$=\psi_{3,2'}=\psi_{4,1'}=\psi_{4,5'}=\psi_{5,1'}=\psi_{5,4'}=0,$$
$$\psi_{2,5'}=\psi_{3,4'}=\psi_{4,3'}=\psi_{5,2'}=1,$$
$$\psi_{2,4'}=\psi_{3,5'}=\psi_{4,2'}=\psi_{5,3'}=2.$$

A \emph{group divisible design} $\mathcal{D}$ with parameters $(n,m;k;\lambda_1,\lambda_2)$, denoted by $GDD(n,m;k;$ $\lambda_1,\lambda_2)$, is an
ordered triple $(\mathcal{P},\mathcal{G},\mathcal{B})$, where $\mathcal{P}$ is a set of points, $\mathcal{G}$ is a partition of $\mathcal{P}$ into $m$ sets of size $n$, each set being called a \emph{group}, and $\mathcal{B}$ is a collection of $k$-sbusets, called blocks of $\mathcal{P}$ so that each pair of points from the same group occurs in exactly $\lambda_1$ blocks and each pair of points from different groups occurs in exactly $\lambda_2$ blocks. The triple $\mathcal{I}=(\mathcal{P},\mathcal{B},I)$ of a group divisible design is an incidence structure with the natural incidence relation $I$, and  the dual incidence structure $\mathcal{I}^*=(\mathcal{B},\mathcal{P},I^*)$. If there exists a partition $\mathcal{G}'$ of $\mathcal{B}$ such that the triple $(\mathcal{B},\mathcal{G}',\mathcal{P})$ is a $GDD(n,m;k;\lambda_1,\lambda_2)$, then we call that the $\mathcal{D}$ is a group divisible design with \emph{dual property} holding parameters $(n,m;k;\lambda_1,\lambda_2)$ and we denote such a design by $GDDDP(n,m;k;\lambda_1,\lambda_2)$.

Let $d \geq  2$ be an integer and let $q$ be a prime power. Let $V$ be a vector space of dimension
$d$ over $GF(q)$, the finite field with $q$ elements. We define the set of non-zero vectors
in $V$ as the point set $\mathcal{P}$ and the set of affine hyperplanes in $V$ as the block set $\mathcal{B}$, i.e.,
$\mathcal{P}= \{x\in  V | x \neq 0\}$ and $ \mathcal{B}= \{x + H | H$ is a hyperplane in $V$ and  $x\in V \setminus  H\}$.

We make a partition $\mathcal{G}$ of $\mathcal{P}$ such that the collinear non-zero vectors in $V$ belong
to the same group in $\mathcal{G}$. Note that each group in $\mathcal{G}$ has size $q-1$. Then $(\mathcal{P}; \mathcal{G}; \mathcal{B})$ is
a $GDD(n,m; k; 0; \lambda_2)$, where $n = q-1, m = \frac{q^d-1}{q-1}$
 is the number of projective points
in $V$, $k = q^{d-1}$ is the number of affine hyperplanes containing a given non-zero vector,
and $\lambda_2 = q^{d-2}$ is the number of affine hyperplanes containing two given non-zero and
non-collinear vectors.

We look at the dual incidence structure $\mathcal{I}^*$ of $\mathcal{I} = (\mathcal{P}; \mathcal{B}; I)$. We
make a partition $\mathcal{G}'$ of $\mathcal{B}$ such that the parallel affine hyperplanes belong to the same
group in $\mathcal{G}'$. Then $(\mathcal{B}; \mathcal{G}';\mathcal{P})$ becomes a $GDD(n;m; k; 0; \lambda_2)$, where $n = q- 1$, $m = \frac{q^d-1}{q-1}$
is the number of $d-1$ dimensional subspaces in $V$, $k = q^{d-1}$ is the number of non-zero
vectors in an affine hyperplane, and $\lambda_2 = q^{d-2}$ is the number of non-zero vectors in the
intersection of two given non-parallel affine hyperplanes.

This shows that $\mathcal{D}(d; q):= (\mathcal{P};\mathcal{G}; \mathcal{B})$ is a $GDDDP(n;m; k; 0; \lambda_2)$, where $n, m, k, \lambda_2$
are given above. We denote $\Gamma(d, q) := \Gamma(\mathcal{D}(d; q))$ as the point-block incidence graph of $\mathcal{D}(d; q)$. It is clear that the general linear group $GL(d, q)$ acts as a group of automorphism of the graph $\Gamma(d, q)$.
By \cite[Proposition 4]{qdk-1},  the point-block incidence graphs $\Gamma(d, q)$ of classical group divisible
designs with the dual property $\mathcal{D}(d; q)$ are 2-arc-transitive dihedrants.

For any $r|q-1$, let $N\leq Z(GL(d,q))\cong \mathbb{Z}_{q-1}$ and $|N|=(q-1)/r$. Let
$\Gamma(d, q, r)$ be the quotient graph of $\Gamma(d,q)$ induced by $N$.

We use  $X'(3, 2)$ to denote the graph with the vertex set
$\{i|i\in \mathbb{Z}_{14}\}\cup \{i'|i'\in \mathbb{Z}_{14}\}$ where $i$ and $j'$ are adjacent if
and only if $j'-i\in \{0,1,9,11\}$.

Let $V(2,2)$ be the 2-dimensional vector space over $\mathbb{F}_2$. Let $\K_{4,4}$
be the complete bipartite graph with the bipartition $V(\K_{4,4})=U\cup W$, where
$U=\{\alpha|\alpha \in V(2,2)\}$
and
$W=\{\alpha'|\alpha \in V(2,2)\}$.
Define the cover $X(2,2)=\K_{4,4}\times_\psi \mathbb{Z}_2$ with the voltage assignment
$\psi$: $\psi(\alpha,\beta')=\alpha \beta^T=a_1b_1+a_2b_2$, for all
$\alpha=(a_1,a_2)$ and $\beta=(b_1,b_2)\in V(2,2)$.


\subsection{Some lemmas}

The first lemma    can be easily proved.

\begin{lemma}\label{2dt-volt-lem1}
Let $X\mapsto Y$ be a regular cyclic covering of a connected graph such that some $2$-distance-transitive group $G\leq \Aut(X)$
projects along $X\mapsto Y$. Then there exists a regular prime cyclic covering
$X'\mapsto Y$ such that some $2$-distance-transitive group $G'\leq \Aut(X')$ projects along
$X'\mapsto Y$.

\end{lemma}

\begin{lemma}{\rm (\cite[Lemma 16.3]{Biggs})}\label{cayley-1}
A graph $\Gamma$ is a Cayley graph on a group $T$ if and only if $\Aut(\Gamma)$ has a regular
subgroup isomorphic to $T$.
\end{lemma}





\begin{theo}{\rm (\cite[Theorem 1.2]{DMM-2008})}\label{2at-dih-du}
Let $n\geq 3$ and let $X$ be a connected $2$-arc-transitive Cayley graph of a dihedral group of order $2n$. Then one of the following occurs:

 \begin{itemize}

\item[(1)] $X$ is a basic graph and is isomorphic to one of the following graphs: $C_{2n}$, $n$ a prime; $\K_{2n}$;
$\K_{n,n}$; $B(H_{11})$ or $B'(H_{11})$; $B(PG(d,q))$ or $B'(PG(d,q))$, where $n=(q^d-1)/(q-1)$, $d\geq 2$, and $q$ is a prime power; or

\item[(2)]   $X$ is not a basic graph and either $X$ is isomorphic to $\K_{n,n}-n\K_2$, or there exists an odd
prime power $q$ such that $n=q+1$ and $X$ is isomorphic to $\K_{q+1}^{2d}$, where $d$ is a divisor of
$\frac{q-1}{2}$ if $q\equiv 1 \pmod 4$, and a divisor of $q-1$ if $q\equiv 3 \pmod 4$, respectively.

\end{itemize}

\end{theo}

\begin{lemma}{\rm (\cite[Lemma 2.5]{Li-abeliancay-2008},\cite[Theorem 4.1]{Praeger-4})}\label{2at-2at-1}
Let  $\Gamma$ be a connected $G$-arc-transitive graph. Suppose that $N$ is a normal subgroup of $G$ such that
$\Gamma$ is an $N$-cover of the normal quotient graph $\Gamma_N$. Then
$\Gamma$ is $(G,2)$-arc-transitive if and only if $\Gamma_N$ is $(G/N,2)$-arc-transitive.
\end{lemma}

The classification of    primitive permutation groups  that contain a
cyclic regular subgroup was independently obtained by Jones \cite{Jones-2002} and Li \cite[Corollary 1.2]{LCH-abelianregular-2003}.
Moreover, by \cite[Theorem 1.2]{LP-circulant-2012}, every   quasiprimitive group with a regular cyclic  subgroup is  primitive.
Hence the family of  quasiprimitive groups with a regular cyclic  subgroup is also  completely  determined in the following lemma.

\begin{theo}{\rm (\cite{Jones-2002},\cite{LCH-abelianregular-2003},\cite{SLZ-2014})}\label{regcyc}
Let $G$ be a quasiprimitive permutation group on $\Omega$ which contains a regular cyclic  subgroup $T$ of degree $n$. Then $G$ is primitive on $\Omega$,
and either $n=p$ is prime and $G\leq AGL(1,p)$, or $G$ is $2$-transitive, listed in Table \ref{cycregular}.

\end{theo}

\begin{table}[]\caption{Quasiprimitive  groups containing regular cyclic subgroups }\label{cycregular}
\medskip
\centering
\begin{tabular}{|c|c|c|c|c|c|}
\hline     $G$   & $G_u$   & $n$ & Condition & $3$-transitive? \\
\hline  $A_n$     &  $A_{n-1}$ & &  $n \geq 5$ is odd   & Yes\\
\hline  $S_n$     &   $S_{n-1}$ & $n\geq 4$ & & Yes\\
\hline   $PGL(2,q).o$ &  $[q]:GL(1,q)$ &  $\frac{q^2-1}{q-1}$  & $o\leq P\Gamma L(2,q)/PGL(2,q)$  & Yes \\
\hline   $PGL(d,q).o$, $d\geq 3$ &  $[q^{d-1}]:GL(d-1,q)$ &  $\frac{q^d-1}{q-1}$ & $o\leq P\Gamma L(d,q)/PGL(d,q)$ & No \\
\hline  $PSL(2,11)$ &   $A_5$ & 11 & & No \\
\hline  $M_{11}$ &   $M_{10}$ & 11& & Yes \\
\hline  $M_{23}$ &   $M_{22}$ & 23& & Yes \\
  \hline
\end{tabular}
\end{table}

The following result about primitive permutation groups  that contain an
element with exactly two equal cycles is due to M\"{u}ller.

\begin{theo}{\rm (\cite[Theorem 3.3]{Muller-2018})}\label{bicirculant-primitive-2}
Let $G$ be a primitive permutation group of degree $2n$ that contains an
element with exactly two cycles of length $n$. Then one of the following holds, where $G_0$ denotes the stabiliser of a point.
\begin{enumerate}[$(1)$]
\item  (Affine action) $\mathbb{Z}_2^m\lhd G\leqslant \AGL(m,2)$ is an affine permutation group, where $n=2^{m-1}$. Further, one of the following holds.
 \begin{enumerate}[{\rm (a)}]
\item $n= 2$, and $ G_0= \GL(2,2)$;
\item $n=2$, and $ G_0= \GL(1,4)$;
\item $n=4$, and $ G_0= \GL(3,2)$;
\item $n=8$, and $G_0\in \{\mathbb{Z}_5:\mathbb{Z}_4,\Gamma L(1,16),(\mathbb{Z}_3\times \mathbb{Z}_3):\mathbb{Z}_4,\Sigma L(2,4),\Gamma L(2,4),A_6,GL(4,2),\\(S_3\times S_3):\mathbb{Z}_2,S_5,S_6,A_7\}.$\footnote{We note that \cite[Theorem 3.3]{Muller-2018} gives $|\Gamma L(1,16):G_0|=3$ as one of the possibilities but there is a unique such group, namely the $\mathbb{Z}_5:\mathbb{Z}_4$ that we list here.}
\end{enumerate}
\item (Almost simple action) $S\leqslant G\leqslant \Aut(S)$ for a  nonabelian simple group $S$, and one of the following holds.
 \begin{enumerate}[{\rm (a)}]
\item  $n\geqslant 3$ and  $A_{2n}\leqslant G\leqslant S_{2n}$ in its natural action;
\item $n=5$ and $A_5\leqslant G\leqslant S_5$ in the action on the set of  $2$-subsets of $\{1,2,3,4,5\}$;
\item $n=(q^d-1)/2(q-1)$  and $\PGL(d,q)\leqslant G\leqslant  P \Gamma L(d,q)$ for some odd  prime power $q$ and $d\geqslant 2$ even;
\item $n=11$ and $M_{22}\leqslant G\leqslant \Aut(M_{22})$;
\item $n=6$ and $G=M_{12}$;
\item $n=12$ and $G=M_{24}$.
\end{enumerate}
\end{enumerate}
\end{theo}

\begin{lemma}{\rm( \cite[Theorem 4.1]{Taylor-1})}\label{2-trans-group}
Suppose that $G$ is a finite $2$-transitive  permutation group on a set of size $n$.

(1) Suppose that $G$ is almost simple and $T\unlhd G\leq Aut(T)$ with $T$ a non-abelian simple group. Then one of the following cases  holds:

(a) alternating: $T=A_n$, with $n\geq 5$;

(b) Lie type: $T=PSL(d,q)$, $d\geq 2$, $n=(q^d-1)/(q-1),(d,q)\neq (2,2),(2,3)$, and
two representations if $d>2$;

(c) Lie type: $T=PSU(3,q)$,  $n=q^3+1,q>2$;

(d) Lie type: $T=^2B_2(q)$,  $n=q^2+1,q=2^{2e+1}>2$;

(e) Lie type: $T= ^2G_2(q)$,  $n=q^3+1,q=3^{2e+1}$;

(f)  $T=Sp(2d,2)$, $d\geq 3$ and   $n=2^{2d-1}\pm 2^{n-1}$;

(g)  $T=PSL(2,11)$,  $n=11$, two representations;

(h)  $T=A_7$,  $n=15$, two representations;

(i)  $T=M_n$,  $n\in \{11,12,22,23,24\}$, two representations for $M_{12}$;

(j)  $T=M_{11}$,  $n=12$;

(k)  $T=HS$ the Higman-Sims group,  $n=176$, two representations;

(l)  $T=Co_3$ the third Conway group,  $n=276$.

(2) $G$ has a regular normal subgroup $T$ which is elementary abelian of order $n=p^d$, where
$p$ is a prime, $G$ may be identified with a group of affine transformations
$x\mapsto x^g+c$ of $F_p^d$, where $g\in G_0$, and one of the following cases holds.

(a) $G\leq A\Gamma L(1,n)$, the group of affine transformations of $F_n$ extended by the automorphisms of $F_n$;

(b) $SL(d,q)\unlhd G_0$, $n=q^d,d>2$;

(c) $Sp(d,q)\unlhd G_0$, $n=q^d$;

(d) $G_2(q)' \unlhd G_0$, $n=q^6,q $ even;

(e) $G_0\cong A_6$ or $A_7$, $n=16$;

(f) $SL(2,3)\unlhd G_0$ or $SL(2,5)\unlhd G_0$ , $n=p^2$, where $p\in \{5,7,11,19,23,29,59\}$, or $n=3^4$;

(g) $G_0$ has a normal extraspecial supbgroup $E$ of order $2^5$, $G/E$ is isomorphic to a subgroup of $S_5$, and $n=3^4$;

(h) $G_0=SL(2,13)$, $n=3^6$.

In each case where $G$ has two representations of degree $n$, these representations are interchanged by an outer automorphism. Moreover, $G_2(q)=G_2(q)'$ when $q\neq 2$ and
$G_2(2)'=PSU(3,3)$ has index $2$ in $G_2(2)$.

\end{lemma}

\begin{propo}{\rm (\cite{Cameron-1})}\label{3trans-infint-1}
Let $G$ be a  $3$-transitive permutation group of degree at least $4$.  Then one of the following occurs.

{\rm (1)} The socle of $G$ is $3$-transitive;

{\rm (2)} $PSL(2,q)\unlhd G\leq P\Gamma L(2,q)$ with natural action on
the projective line $PG(1,q)$ of degree $q+1$, for odd $q\geq 5$, and with the socle of $G$ being isomorphic to $PSL(2,q)$ and not $3$-transitive;

{\rm (3)} $G=AGL(m,2)$, $m\geq 3$;

{\rm (4)} $G=\mathbb{Z}_2^4:A_7<AGL(4,2)$;

{\rm (5)} $G=PGL(2,3)\cong AGL(2,2)\cong S_4$ of degree $4$;

\end{propo}

\section{ Proof of main Theorem }

In this section, we will prove our main theorem by a series of lemmas.



The first lemma provides some  properties   of generalized quaternion groups and these facts will be used in the following argument.
\begin{lemma}\label{lem-gq-norm-1}
Let $Q_{4n}=\langle a,b| a^{2n}=1,b^2=a^{n}, b^{-1}ab=a^{-1}\rangle$ where $n\geq 2$. Then $Q_{4n}$ is not abelian, and the following statements hold.

 \begin{itemize}

\item[(1)] Every element of $Q_{4n}$ can be written as $a^ib^j$ where $1\leq i<2n $ and $j=0,1$.
Moreover, $(a^kb)^2=a^n$, $a^kba^m=a^{k-m}b$, and $a^kba^mb=a^{k-m+n}$.

\item[(2)]   $Q_{4n}$ has a unique element of order $2$, that is, the element $b^2=a^n$.

\item[(3)]  $Z(Q_{4n})=\langle a^n\rangle \cong \mathbb{Z}_2$, and
$Q_{4n}/Z(Q_{4n})\cong D_{2n}$.

\item[(4)]  Assume  $n$ is odd. Then   all the normal subgroups of $Q_{4n}$ are the following:
$1$, $Q_{4n}$, $\langle a^{2n/k}\rangle$, where $k$ is a divisor of $2n$.
Moreover,  the only index $2$  normal subgroup of $Q_{4n}$ is
 $\langle a\rangle$.

\item[(5)] Assume  $n$ is even. Then  all the normal subgroups of $Q_{4n}$ are the following:
$1$, $Q_{4n}$, $\langle a^{2n/k}\rangle$, where $k$ is a divisor of $2n$,  $\langle a^{2},b\rangle$
and $\langle a^{2},ab\rangle$. Moreover,
  $Q_{4n}$ has three normal subgroups of index $2$:
 $\langle a\rangle$,  $\langle a^{2},b\rangle$
and $\langle a^{2},ab\rangle$.

\item[(6)]   If $n>2$ is even, then    $\langle a^{2},b\rangle$
and $\langle a^{2},ab\rangle$ are generalized quaternion groups of order $2n$.

\end{itemize}
\end{lemma}



We give a remark on the connected Cayley  graphs over  generalized quaternion groups.

\begin{rem}\label{rem-1}
{\rm Let $T:=\langle a,b|a^{2n}=1,a^n=b^2,a^b=a^{-1}\rangle \cong Q_{4n}$, where $n\geq 2$. Let $\Gamma=\Cay(T,S)$ be a connected Cayley graph. Then $T=\langle S\rangle$ and $S=S^{-1}$. If
$|S|=2$, then as $T$ has a unique element of order $2$,
$S=\{x,x^{-1}\}$ where $x$ is an element of order at least 3. However, in this case
 $T$ is a cyclic group, a   contradiction. If
$|S|=3$, then $S=\{a^n,x,x^{-1}\}$ where $x$ is an element of order at least 3.
If $x\in \langle a\rangle$, then   $\langle S\rangle$ is cyclic, and  $T\neq \langle S\rangle$.
If $x=a^ib$, then $x^2=a^n$, and so $\langle S\rangle=\langle x\rangle$ is cyclic,
 and so
$\Gamma$ is disconnected. Therefore  $|S|\geq 4$.

}
\end{rem}

Note that  the cyclic subgroup  $\langle a\rangle\cong \mathbb{Z}_{2n}$ has 2 orbits on the vertex set of $\Gamma$.  By \cite[Lemma  3.1]{DGJ-2019}, we have the following lemma.

\begin{lemma}\label{quo-1}
Let $\Gamma$ be a connected $G$-arc-transitive Cayley  graph over the group    $T=\langle a,b|a^{2n}=1,a^n=b^2,a^b=a^{-1}\rangle$, where $n\geq 2$.
Let $H=\langle a\rangle\cong \mathbb{Z}_{2n}$.
Let $H_0$ and $H_1$ be the two orbits of $H$ on $V(\Gamma)$ and
let $\mathcal{B}$ be a $G$-invariant partition of $V(\Gamma)$. Then the following statements hold.

\begin{enumerate}[$(1)$]
\item Either all elements of $\mathcal{B}$ are subsets of $H_0$ or $H_1$;
or  $B\cap H_0\neq \varnothing$ and $B\cap H_1\neq \varnothing$ for every   $B\in \mathcal{B}$.

\item If $B\in \mathcal{B}$ and  $B\cap H_i\neq \varnothing$ for some $i$, then $B\cap H_i$ is a block for $H$ on $H_i$.
\end{enumerate}

\end{lemma}

We give a  reduction result.

\begin{theo}\label{bic-reduction-th2}
Let $\Gamma$ be a connected $(G,2)$-distance-transitive Cayley  graph over the generalized quaternion group    $T=\langle a,b|a^{2n}=1,a^n=b^2,a^b=a^{-1}\rangle$, where $n\geq 2$
and $T\leq G\leq \Aut(\Gamma)$.
 Suppose that $\Gamma\ncong \K_{m[f]}$ for any  $m\geq 3$ and $f\geq 2$.
 Let   $N$ be a  normal subgroup of $G$ maximal with respect to having at least $3$ orbits. Then
$\Gamma$ is a cover of $\Gamma_N$ which
 is a $(G/N,s)$-distance-transitive Cayley graph over $T/N$, where $s=\min\{2,\diam(\Gamma_N)\}$, and $G/N$ is faithful and either quasiprimitive or bi-quasiprimitive on $V(\Gamma_N)$.

Moreover,   $N=\langle a^i\rangle$ where $i$ is a divisor of $2n$ and $i\neq 1,2n$,  and one of the following holds:

\begin{enumerate}[$(1)$]
\item    $T/N\cong \mathbb{Z}_2\times \mathbb{Z}_2$, $i=2$ and $n$ is even;

\item    $T/N\cong \mathbb{Z}_4$, $i=2$ and $n$ is odd;

\item      $T/N\cong Q_{2i}$, $i\geq 3$ is even  and $i$ is not a divisor of $n$;

\item      $T/N\cong D_{2i}$, $i\geq 3$  and $i$ is  a divisor of $n$.
\end{enumerate}

\end{theo}
\proof
Since $\Gamma$ is a $(G,2)$-distance-transitive graph, it  is also locally $(G,2)$-distance-transitive, and hence \cite[Lemma 5.3]{DGLP-locdt-2012} applies.
Due to
$N$ is intransitive on $V(\Gamma)$ and using the  $G$-arc-transitivity of $\Gamma$, we know that each non-trivial $N$-orbit does
not contain any edge of $\Gamma$.
As     $N$ has at least 3 orbits in $V(\Gamma)$,  the fact $\Gamma\ncong \K_{x[y]}$ for any   $x\geq 3$ and $y\geq 2$ implies  that only \cite[Lemma 5.3]{DGLP-locdt-2012}  (iv) occurs. Hence $N$ is
semiregular on the vertex set  and  $\Gamma$ is a cover of $\Gamma_N$.

Since $\Gamma$ is    $G$-arc-transitive,  we can easily show that $\Gamma_N$ is $G/N$-arc-transitive. Assume that  $\Gamma_N$ is a non-complete graph.
Let $(C_1,C_3)$ and $(C_1',C_3')$ be two pairs of vertices of  $\Gamma_N$ such that $d_{\Gamma_N}(C_1,C_3)=d_{\Gamma_N}(C_1',C_3')=2$. Then there exist $c_i\in C_i$ and $c_i'\in C_i'$ such that
$(c_1,c_3)$ and $(c_1',c_3')$ are two pairs of vertices  of $\Gamma$ with $d_{\Gamma}(c_1,c_3)=d_{\Gamma}(c_1',c_3')=2$. Since $\Gamma$ is $(G,2)$-distance-transitive, there exists $\alpha \in G$
such that   $(c_1,c_3)^\alpha=(c_1',c_3')$. Hence $(C_1,C_3)^\alpha=(C_1',C_3')$. In particular, $\alpha $ induces an element of $ G/N$ that maps $(C_1,C_3)$ to $(C_1',C_3')$.
Therefore $\Gamma_N$ is   $(G/N,2)$-distance-transitive.

Again using the fact that   $N$ is a  normal subgroup of $G$ maximal with respect to having at least $3$ orbits,
all normal subgroups of the quotient group $G/N$ are transitive or have exactly two orbits on $V(\Gamma_N)$. Thus $G/N$ is quasiprimitive or bi-quasiprimitive on $V(\Gamma_N)$.

Since $T$ acts  regularly on $V(\Gamma)$ and  $H:=\langle a\rangle\cong \mathbb{Z}_{2n}$ is an index 2 subgroup of $T$, it follows that $H$ has 2 orbits on $V(\Gamma)$.
Let $H_0$ and $H_1$ be the two $H$-orbits on $V(\Gamma)$.

Let $\mathcal{B}=\{B_1,\ldots,B_t\}$ be the set of $N$-orbits where  $t\geqslant 3$.
By Lemma \ref{quo-1},  either

\begin{enumerate}[{\rm (i)}]
\item  all elements of $\mathcal{B}$ are subsets of $H_0$ or $H_1$ and $N\leqslant H$ is a cyclic group; or
\item  $B\cap H_0\neq \varnothing$ and $B\cap H_1\neq \varnothing$ for each   $B\in \mathcal{B}$, and  $N$       has a cyclic index $2$ normal subgroup $H\cap N$.

\end{enumerate}

Suppose  that (ii) holds. Then  $N$       has a cyclic index $2$ normal subgroup $H\cap N$, and  hence $N\nleq H$.
Since $N$ is semiregular on $V(\Gamma)=T$, we have $N\leq T$, and so $N$ is a normal subgroup of $T$ with at least 3 orbits.

If   $n$ is odd,  then by  Lemma \ref{lem-gq-norm-1},    all the normal subgroups of $T$ are the following:
$1$, $T=Q_{4n}$, $\langle a^{2n/k}\rangle$, where $k$ is a divisor of $2n$, contradicting the fact that $N$ is a normal subgroup of $T$.

We consider that      $n$ is even.   Then    all the normal subgroups of $T$ are the following:
$1$, $T=Q_{4n}$, $\langle a^{2n/k}\rangle$, where $k$ is a divisor of $2n$,  $\langle a^{2},b\rangle$
and $\langle a^{2},ab\rangle$.
Since $N\nleq H$ and $N\neq T$,  it follows that $N=\langle a^{2},b\rangle$
or   $\langle a^{2},ab\rangle$.
Due to     $n>2$ is even,     both $\langle a^{2},b\rangle$
and $\langle a^{2},ab\rangle$ are generalized quaternion groups of order $2n$.
As a result   $|T:N|=2$, and $N$ has 2 orbits on $V(\Gamma)$, again a contradiction.
Hence case (ii) does not occur.

Now assume  that (i) holds. Then $N\leq H$.  Since $N$ has at least 3 orbits on $V(\Gamma)$, we have $N<H$.
Consequently    $N=\langle a^i\rangle$ where $i$ is a divisor of $2n$ and $i\neq 1,2n$, and $|N|=2n/i$.
Since $H$ is cyclic, $N$ is a characteristic subgroup of $H$, and so $N$ is an normal subgroup of $T$.
Moreover, $T/N$ is regular on $V(\Gamma_N)$, and so by Lemma \ref{cayley-1},  $\Gamma_N$ is a Cayley graph over the group $T/N$.

Since $N$ is an normal subgroup of $T$, it follows from    Lemma \ref{lem-gq-norm-1} that   if $i=2$, then $|T/N|=4$, and  $T/N\cong \mathbb{Z}_2\times \mathbb{Z}_2$ when $n$ is even,
and $T/N\cong \mathbb{Z}_4$ when   $n$ is  odd;
for $i\geq 3$, if $i$ is not a divisor of $n$, then $i$ is even, and $T/N\cong Q_{2i}$; and if
 $i$ is  a divisor of $n$, then  $T/N\cong D_{2i}$.

We conclude the proof.
\qed


\begin{lemma}\label{2dt-cay-quasi-th1}
Let $\Gamma$ be a connected $G$-arc-transitive  Cayley  graph over the group    $T=\langle a,b|a^{2n}=1,a^n=b^2,a^b=a^{-1}\rangle$, where $n\geq 2$
and $T\leq G\leq \Aut(\Gamma)$.
Suppose that  $G$ is  quasiprimitive on  $V(\Gamma)$.
Then $G$ is  $2$-transitive on $V(\Gamma)$ and  $\Gamma$ is the complete graph $\K_{4n}$.

\end{lemma}
\proof
Assume  that $G$ is quasiprimitive on $V(\Gamma)$. If $G$ acts  primitively on $V(\Gamma)$, then since $T$ is a generalized quaternion group, it follows from  \cite{LPX-2021} and \cite{Scott-1957} that
$G$ is  2-transitive on $V(\Gamma)$.
Hence $\Gamma$ is a complete graph, and so it is isomorphic to $\K_{4n}$.

Assume that $G$ is not  primitive on $V(\Gamma_N)$.
Let $\mathcal{B}=\{B_1,\ldots,B_t\}$ be the set of $N$-orbits where  $t\geqslant 3$.
Let $H:=\langle a\rangle\cong \mathbb{Z}_{2n}$ which is an index 2 subgroup of $T$. Then   $H$ has 2 orbits on $V(\Gamma)$.
Let $H_0$ and $H_1$ be the two $H$-orbits on $V(\Gamma)$.

Since $G$ is quasiprimitive on $\Omega$, it follows that $G$ acts faithfully on $\mathcal{B}$. Moreover, the maximality of $\mathcal{B}$
implies that $G$ is primitive on $\mathcal{B}$.
By Theorem \ref{bic-reduction-th2},   for each $B\in \mathcal{B}$,  either $B\subseteq H_0$ or $B\subseteq H_1$.
Then $|B|$ divides $|H|$, without loss of generality,  we
assume that   $H_0=B_1\cup \cdots\cup B_r$ and
$H_1=B_{r+1}\cup \cdots\cup B_{2r}$.
As  $H_0$ and $H_1$ are the two orbits of $H$, it follows that  $H$ is transitive on both sets $\{B_1,\ldots,B_r\}$ and $\{B_{r+1},\ldots,B_{2r}\}$.
Since $G$ is faithful on $\mathcal{B}$, $H$ is
faithful  on $\mathcal{B}$.
Further, due to $H$ is a cyclic group, it is regular on each orbit on $\mathcal{B}$.
Thus $|\mathcal{B}|=2|H|=|\Omega|$,
contradicting the fact that $\mathcal{B}$ is a non-trivial  $G$-invariant partition of $\Omega$.
\qed

Let $\Gamma$ be a  $2$-distance-transitive graph.
If  $\Gamma$ has girth  $3$, then
 $\Gamma$ has two types of 2-arcs, and it  is  not   $2$-arc-transitive.
If  $\Gamma$ has girth at least $5$, then for any two vertices $u,v$ with distance 2, there is a unique 2-arc between $u$ and $v$, and so
 $\Gamma$ is $2$-arc-transitive.
For the case that    $\Gamma$ has girth  $4$,
 $\Gamma$ may be and may be  not   $2$-arc-transitive.
The following lemma shows that if   $\Gamma\cong \K_{2n,2n}-2n\K_2$  and  $G$ is  bi-quasiprimitive on  $V(\Gamma)$,
then $\Gamma$ is  $(G,2)$-arc-transitive.

Since $\K_{2n,2n}-2n\K_2$ has valency $2n-1$, if
$ \K_{2n,2n}-2n\K_2$ is a connected Cayley  graph over an order $4n$ generalized quaternion group, then by Remark \ref{rem-1},  $n\geq 3$.

\begin{lemma}\label{2dt-knn-k2}
Suppose that   $\Gamma\cong \K_{2n,2n}-2n\K_2$ is a $(G,2)$-distance-transitive Cayley  graph over the group    $T=\langle a,b|a^{2n}=1,a^n=b^2,a^b=a^{-1}\rangle$, where $n\geq 3$
and $T\leq G\leq \Aut(\Gamma)$.
If  $G$ is  bi-quasiprimitive on  $V(\Gamma)$,
then $\Gamma$ is  $(G,2)$-arc-transitive.

\end{lemma}
\proof  Suppose that  $G$ is bi-quasiprimitive on $V(\Gamma)$.
Then $G$ has a minimal normal subgroup $M$ that has exactly two orbits on $V(\Gamma)$.
Since  $\Gamma$ is $G$-arc-transitive and connected,
each  $M$-orbit does not  contain any edge of $\Gamma$.
Thus   the two $M$-orbits form the two bipartite halves of $\Gamma$.
In particular,   all intransitive  normal subgroups of $G$ have the same orbits.
Let $\Delta_0$ and $\Delta_1$ be the two bipartite halves of $\Gamma$.
Let $G^+=G_{\Delta_0}=G_{\Delta_1}$.
Since $G$ is transitive on the vertex set of $\Gamma$, it follows that  $G=\langle G^+, \sigma \rangle$	 for some element $\sigma\in G\setminus G^+$ with $\sigma^2\in G^+$.

Suppose that     $G^+$ acts unfaithfully on each $\Delta_i$.
Let
$K_0$ be the kernel of $G^+$ on $\Delta_0$ and $K_1$ be the kernel of $G^+$ on $\Delta_1$.
Since   $G^+$ acts unfaithfully on each $\Delta_i$, it follows that   $K_i$ is not trivial for $i=0,1$.   As   $\sigma\in G$  swaps   $\Delta_0$ and $\Delta_1$,  $K_i^\sigma=K_{1-i}$  and hence
$1 \neq K_0\times K_1$ is a normal subgroup of  $G$.
Due to  $G$ is bi-quasiprimitive on $V(\Gamma)$ and $K_0\times K_1$ is not transitive, it follows that
$K_0\times K_1$ has two orbits on $V(\Gamma)$, namely $\Delta_0$ and $\Delta_1$. As a result   $K_0$  fixes point-wise   $\Delta_0$ and is transitive on $\Delta_1$. Since $\Gamma$ is $G$-arc-transitive it follows that $\Gamma$ is a complete bipartite graph, which is a contradiction.

Thus    $G^+$ acts faithfully on each $\Delta_i$.
Since $\Gamma$ is a Cayley graph over the group $T$, it follows that  $T$ is regular on $V(\Gamma)$. Due to    $H=\langle a\rangle\cong \mathbb{Z}_{2n}$ is an index 2 subgroup of $T$, it follows that $H$ has 2 orbits on $V(\Gamma)$.
Let $H_0$ and $H_1$ be the two $H$-orbits on $V(\Gamma)$.

Suppose first that  the two orbits of the cyclic subgroup $H$ are the two bipartite halves of $\Gamma$.

Assume that $G^+$ acts  imprimitively on $\Delta_0$. Let $\mathcal{B}_0$ be  a maximal $G^+$-invariant partition  of $\Delta_0$. Then $\mathcal{B}_0^\sigma$ is a maximal $G^+$-invariant partition  of $\Delta_1$.
Let $\mathcal{C}=\mathcal{B}_0 \cup \mathcal{B}_0^\sigma$. Then $G$ leaves invariant this partition
$\mathcal{C}$ of the vertex set, as $\mathcal{B}_0^{\sigma^2}=\mathcal{B}_0$. Since  $G$ is bi-quasiprimitive on $V(\Gamma)$ and
$|\mathcal{C}|\geqslant 4$, it follows that $G$  acts faithfully on $\mathcal{C}$.
Let $M_0$ be the kernel of $H$ acting on $\mathcal{B}_0$ and $M_1$ be the kernel of $H$ acting on $\mathcal{B}_0^\sigma$.
Since $H$ acts transitively on $\mathcal{B}_0$ and  $\mathcal{B}_0^\sigma$, it follows that  $|M_0|= |H|/|\mathcal{B}_0|=|H|/|\mathcal{B}_0^{\sigma}|=|M_1|$. Then due to  $H$ is a cyclic group and has a unique subgroup of each order, it follows that
$M_0=M_1$, that is, the kernel of $H$ on $\mathcal{B}_0$ is the same as the kernel on $\mathcal{B}_0^\sigma$, and is hence in the kernel of $G$ on $\mathcal{C}$. However, $G$ is faithful on $\mathcal{C}$, and so $H$ acts faithfully on $\mathcal{B}_0$. Thus $|\mathcal{B}_0|=|H|$, contradicting $G^+$ being imprimitive on $H_0$.

Therefore  $G^+$ is primitive on each $\Delta_i$.
Since $G^+$ contains the cyclic subgroup $H$ that is transitive on each $G^+$-orbit, it follows from Theorem \ref{regcyc} that
either $|H|=p$ is a prime or $G^+$ is $2$-transitive on $\Delta_i$.
Since   $|H|=2n$,   $G^+$ must be   $2$-transitive on $\Delta_i$
with $|\Delta_i|$ not a prime.  Hence $\Gamma$ is  $(G,2)$-arc-transitive.


Assume  now that the  two orbits of $H$   are not the bipartite halves of $\Gamma$.

Note that $|V(\Gamma)|\geq 4$. If $|V(\Gamma)|= 4$, then  $\Gamma$ is
$ \K_{2,2}$, a contradiction.
Thus  $|V(\Gamma)|>4$. Then $\{\Delta_0,\Delta_1\}$ is the unique $G$-invariant partition of  $V(\Gamma)$ into two equal sized parts.  Since the $\Delta_i$ are not $H$-orbits and $H$ has two orbits of size $2n$, it follows that $H^+:=H\cap G^+$ has index two in $H$ and has two equal sized orbits on each $\Delta_i$. Thus  $G=\langle G^+,H\rangle$.

Suppose first that $G^+$ is primitive on each $\Delta_i$. Since $H^+$ is a cyclic subgroup with exactly two orbits of size $2n$ on $\Delta_i$, the possibilities for $G^+$ are given by Theorem \ref{bicirculant-primitive-2}. In particular, either  $G^+$ is almost simple or $2n=2^{m-1}$, $\mathbb{Z}_2^m \lhd G^+\leqslant \AGL(m,2)$ with $2\leqslant m\leqslant 4$.

If $2n=2^{m-1}$, $\mathbb{Z}_2^m \lhd G^+\leqslant \AGL(m,2)$ with $2\leqslant m\leqslant 4$, then as $n\geq 2$, we have $m\geq 3$. Hence $m=3$ or 4.
Moreover,
$G^+$ has a normal subgroup $N\cong \mathbb{Z}_2^m$ that is regular on each $\Delta_i$.  Now $N$ is characteristic in $G^+$ and so is normal in $G$.  Moreover, since $H$ is cyclic, either $H\cap N=1$ or $H\cap N=\mathbb{Z}_2$.
Due to $N=C_{G^+}(N)$  and $H$ is cyclic, either $|H|=4$ and $C_H(N)=H$, or $C_H(N)= H^+\cap N$.  As  $m= 3$ or 4, it follows that   $|H|=2^m\geqslant 8$ and  $H^+\cap N\neq H^+$,  consequently  $C_H(N)=H^+\cap N$. Thus $H/(H^+\cap N)$ is isomorphic to a  cyclic  subgroup of $\GL(m,2)$ of order $2^m$ or $2^{m-1}$. Since $m=3$ or $4$, it follows that $m=3$, $G^+=\AGL(3,2)$ and $H^+\cap N=2$. As $|H|=8$ and $G^+$ does not have an element of order 8, we have $G=\Aut(\AGL(3,2))$. In particular, given $u\in \Delta_0$ we have that $G_u^+$ is transitive on $\Delta_1$ and so $\Gamma=\K_{8,8}$, a contradiction.

Next suppose that $G^+$ is almost simple. Then either $G^+$ is $2$-transitive on each $\Delta_i$, or $2n=10$ and $A_5\leqslant G^+\leqslant S_5$. In the first case,  $\Gamma$ is $(G,2)$-arc-transitive.

Suppose instead that $2n=10$ and $A_5\leqslant G^+\leqslant S_5$. Then by \cite[p.75]{DM-1}, the action of $G^+$ on $\Delta_i$ is $S_5$ or $A_5$ acting naturally on the set of unordered pairs of $\{1,2,3,4,5\}$.
Let $u\in \Delta_0$.
Then  $G_{u}^+$ has only orbits of size 1, 3 or 6 on both $\Delta_0$ and $\Delta_1$.
Since $G_u^+$ is transitive on $\Gamma(u)$,
it follows that $|\Gamma(u)|=3$ or 6, a contradiction.

We conclude the proof.
\qed

\begin{lemma}\label{2dt-cay-biquasi-th1}
Let $\Gamma$ be a connected $(G,2)$-distance-transitive Cayley  graph over the group    $T=\langle a,b|a^{2n}=1,a^n=b^2,a^b=a^{-1}\rangle$, where $n\geq 2$
and $T\leq G\leq \Aut(\Gamma)$.
If   $G$ is  bi-quasiprimitive on  $V(\Gamma)$,
then $\Gamma$ is one of the following graphs:
$\K_{2n,2n}$,  $ \K_{2n,2n}-2n\K_2$ (with $n\geq 3$),  $B(\PG(d-1,q))$ and $B'(\PG(d-1,q))$, where $d\geq 3$, $q$ is a prime power and $\frac{q^d-1}{q-1}=2n$.
Moreover,
$\Gamma$ is $(G,2)$-arc-transitive.

\end{lemma}
\proof Suppose that  $G$ is bi-quasiprimitive on $V(\Gamma)$.
Then $G$ has a minimal normal subgroup $M$ that has exactly two orbits on $V(\Gamma)$.
Since  $\Gamma$ is $G$-arc-transitive and connected,
it follows that each $M$-orbit does not  contain any   edge of $\Gamma$.
Thus $\Gamma$ is a bipartite graph, and  the two $M$-orbits form the two bipartite halves of $\Gamma$.
Let $\Delta_0$ and $\Delta_1$ be the two bipartite halves of $\Gamma$.

Since $\Gamma$ is a Cayley graph over the group $T$, it follows that $T$ acts regularly on $V(\Gamma)$. As  $H=\langle a\rangle\cong \mathbb{Z}_{2n}$ is an index 2 subgroup of $T$, the group  $H$ has 2 orbits on $V(\Gamma)$.
Let $H_0$ and $H_1$ be the two $H$-orbits on $V(\Gamma)$.

If   the two orbits $H_0$, $H_1$ of the cyclic subgroup $H$ are the two bipartite halves $\Delta_0$, $\Delta_1$ of $\Gamma$, then $\Gamma$ is one of the  graphs
in  \cite[Proposition 5.1]{DGJ-2019}; and if  $H_0$ and  $H_1$ are not  the two bipartite halves $\Delta_0$ and  $\Delta_1$,
then $\Gamma$ is listed in    \cite[Proposition 5.2]{DGJ-2019}.
Thus  $\Gamma$ is one of the following graphs:
$\K_{2n,2n}$,  $ \K_{2n,2n}-2n\K_2$, $B(H(11))$,  $B'(H(11))$, $G(2p,r)$, $B(\PG(d-1,q))$ and $B'(\PG(d-1,q))$, where $d\geq 3$, $q$ is a prime power and $\frac{q^d-1}{q-1}=2n$.

Since $|V(\Gamma)|=|T|$ is divisible by 4, it follows that
the graphs $B(H(11))$,  $B'(H(11))$ and  $G(2p,r)$ do not occur, and so
$\Gamma$ is listed in   the following:
$\K_{2n,2n}$,  $ \K_{2n,2n}-2n\K_2$,  $B(\PG(d-1,q))$ and $B'(\PG(d-1,q))$, where $d\geq 3$, $q$ is a prime power and $\frac{q^d-1}{q-1}=2n$.
Moreover, since $\K_{2n,2n}-2n\K_2$ has valency $2n-1$, it follows from  Remark \ref{rem-1} that $2n-1\geq 4$, that is    $n\geq 3$.

 Lemma \ref{2dt-knn-k2} says that  $ \K_{2n,2n}-2n\K_2$ is $(G,2)$-arc-transitive.
By \cite[Lemma 3.3]{CJS-2dt} and \cite[Lemma 2.11]{HFZY-2025},  the complete bipartite graph $\K_{2n,2n}$  is   $(G,2)$-distance-transitive indicating that it is  also $(G,2)$-arc-transitive.
If $\Gamma$ is  $B(\PG(d-1,q))$ or $B'(\PG(d-1,q))$, where $d\geq 3$, $q$ is a prime power and $\frac{q^d-1}{q-1}=2n$, then
$PGL(d,q)\leq G^+\leq P\Gamma L(d,q)$, and $\Gamma$ is  $(G,2)$-arc-transitive.
\qed

\begin{lemma}\label{2dt-cay-qd-lem-1}
Let $\Gamma$ be a connected $(G,2)$-distance-transitive Cayley  graph over a generalized quaternion group    $T$, where  $T\leq G\leq \Aut(\Gamma)$.
 Suppose that $\Gamma\ncong \K_{m[f]}$ for any  $m\geq 3$ and $f\geq 2$.
Let  $N$ be a  normal subgroup of $G$ maximal with respect to having at least $3$ orbits.
If $\Gamma_N$ is a Cayley graph over the group $T/N\cong D_{2i}$ where $i\geq 3$, then $\Gamma$ is one of the following graphs:
 \begin{itemize}

\item[(1)]        $\K_{2n,2n}-2n\K_2$ with $n\geq 3$;

\item[(2)]        $ X_1(4,q)$ where $q\equiv 3\pmod{4}$ and  $q=n-1$;

\item[(3)]  $\Gamma(d, q, r)$ where $r|q-1$;

\item[(4)] $X(2,2)$;

\item[(5)] $X'(3,2)$.

\end{itemize}

\end{lemma}
\proof Suppose that $\Gamma_N$ is a Cayley graph over the group $T/N\cong D_{2i}$, where $i\geq 3$. Since   $\Gamma\ncong \K_{m[f]}$ for any  $m\geq 3$ and $f\geq 2$, and   $N$ is   a  normal subgroup of $G$ maximal with respect to having at least $3$ orbits, it follows from   Theorem \ref{bic-reduction-th2} that
$\Gamma$ is a cover of $\Gamma_N$ which
 is $(G/N,s)$-distance-transitive, where $s=\min\{2,\diam(\Gamma_N)\}$, and $G/N$ is faithful and either quasiprimitive or bi-quasiprimitive on $V(\Gamma_N)$.

Since $\Gamma_N$ is a  Cayley graph over the group $D_{2i}$,  where $i\geq 3$, it follows from
\cite{HFZY-2025} and \cite{JT-2022}  that
$\Gamma_N$ is   $(G/N,2)$-arc-transitive.
As   $G/N$ acts quasiprimitively or bi-quasiprimitively on $V(\Gamma_N)$,
$\Gamma_N$ is a basic graph, and applying      Theorem \ref{2at-dih-du},
$\Gamma_N$ is  isomorphic to one of the following graphs:
 \begin{itemize}
\item[(1)] $\K_{2i}$;

\item[(2)]   $\K_{i,i}$;

\item[(3)]  $B(PG(d,q))$ or $B'(PG(d,q))$, where $i=(q^d-1)/(q-1)$, $d\geq 2$, and $q$ is a prime power.

\end{itemize}

As $\Gamma_N$ is   $(G/N,2)$-arc-transitive, it follows from   Lemma \ref{2at-2at-1} that
$\Gamma$ is   $(G,2)$-arc-transitive.
If $\Gamma_N$ is isomorphic to $\K_{2i}$, then by \cite[Lemma 4.1]{WJ-2023jctb},
$\Gamma$ is either the graph $\K_{2n,2n}-2n\K_2$ with $n\geq 3$ or $ X_1(4,q)$ where $q\equiv 3\pmod{4}$ and  $q=2i-1$.
Note that $ X_1(4,q)$ has $4(q+1)$ vertices. Hence $4(q+1)=4n$, and $q=n-1$.

For the case that   $\Gamma_N$ is isomorphic to $\K_{i,i}$,  \cite[Lemma 4.3]{WJ-2023jctb} says that
$i=4$ and  $\Gamma$ is the graph $X(2,2)$.

By \cite[Lemma 4.6]{WJ-2023jctb} (1),  $B(\PG(d-1,q))$  does not have  $(G,2)$-arc-transitive  regular cyclic   covers as bicirculants, where $d\geq 3$ and  $q$ is a prime power. Hence  $\Gamma_N$ is not isomorphic to $B(\PG(d-1,q))$.

For the case that   $\Gamma_N$ is isomorphic to $B'(\PG(d-1,q))$, where $d\geq 3$, $q$ is a prime power and $\frac{q^d-1}{q-1}=i$,   \cite[Lemma 4.6]{WJ-2023jctb} (2) indicates  that
$\Gamma$ is either $X'(3,2)$ or the graph $\Gamma(d, q, r)$ where $r|q-1$.
\qed

\begin{lemma}\label{2dt-cay-qd-lem-2}
Let $\Gamma$ be a connected $(G,2)$-distance-transitive Cayley  graph over the group    $T=\langle a,b|a^{2n}=1,a^n=b^2,a^b=a^{-1}\rangle$, where $n\geq 2$
and $T\leq G\leq \Aut(\Gamma)$.
 Suppose that $\Gamma\ncong \K_{m[f]}$ for any  $m\geq 3$ and $f\geq 2$.
Let  $N$ be a  normal subgroup of $G$ maximal with respect to having at least $3$ orbits.
Suppose that $\Gamma_N$ is a Cayley graph over the group  $T/N\cong Q_{2i}$, $i\geq 3$ is even.
If $G/N$ is  bi-quasiprimitive on $V(\Gamma_N)$, then $\Gamma$ is one of the following graphs:
 \begin{itemize}


\item[(1)]   $\K_{q+1}^{2d}$, where   $d\geq 2$   dividing $q-1$,  $q$ is an odd prime power and $2n=d(q+1)$;

\item[(2)] $\Gamma(d, q, r)$ where $r|q-1$;



\item[(4)]  $ X_2(3)$;

\item[(5)] $X(2,2)$;

\item[(6)] $X'(3,2)$.

\end{itemize}

\end{lemma}
\proof    Since   $\Gamma\ncong \K_{m[f]}$ for any  $m\geq 3$ and $f\geq 2$, and   $N$ is   a  normal subgroup of $G$ maximal with respect to having at least $3$ orbits, it follows from   Theorem \ref{bic-reduction-th2} that
$\Gamma$ is a cover of $\Gamma_N$ which
 is $(G/N,s)$-distance-transitive, where $s=\min\{2,\diam(\Gamma_N)\}$.

Suppose that $G/N$ is  bi-quasiprimitive on $V(\Gamma_N)$. It follows from   Lemma \ref{2dt-cay-biquasi-th1} that
$\Gamma_N$ is   $(G/N,2)$-arc-transitive, and $\Gamma_N$ is one of the following graphs:
$\K_{i,i}$,  $ \K_{i,i}-i\K_2$,  $B(\PG(d-1,q))$ and $B'(\PG(d-1,q))$, where $d\geq 3$, $q$ is a prime power and $\frac{q^d-1}{q-1}=i$.
Moreover,  by Lemma \ref{2at-2at-1},
$\Gamma$ is   $(G,2)$-arc-transitive.

If $\Gamma_N$ is isomorphic to $\K_{i,i}$, then by \cite[Lemma 4.3]{WJ-2023jctb},
$i=4$ and  $\Gamma$ is the graph $X(2,2)$.

Assume that  $\Gamma_N$ is isomorphic to $\K_{i,i}-i\K_2$ and  set $|N|=d$.
By Remark \ref{rem-1},   $\Gamma$ has valency at least 4. Since $\Gamma_N$ has the same valency as $\Gamma$, $\Gamma_N$ has valency at least 4, and so $i\geq 5$.
Then applying \cite[Lemma 4.7]{WJ-2023jctb},
 either

\begin{itemize}



\item[(1)]  $i=5$,   $d=3$ and  $ \Gamma=X_2(3)$; or

\item[(2)] $i=q+1\geq 6$ for some odd prime power $q$, and  $\Gamma=\K_{q+1}^{2d}$ for some $d\geq 2$   dividing $q-1$ and $2n=d(q+1)$.
\end{itemize}

By \cite[Lemma 4.6]{WJ-2023jctb} (1),   $\Gamma_N$ is not isomorphic to $B(\PG(d-1,q))$, where $d\geq 3$ and  $q$ is a prime power.

Finally, if    $\Gamma_N$ is isomorphic to $B'(\PG(d-1,q))$, where $d\geq 3$, $q$ is a prime power and $\frac{q^d-1}{q-1}=i$,  then by \cite[Lemma 4.6]{WJ-2023jctb} (2),
$\Gamma$ is either $X'(3,2)$ or the graph $\Gamma(d, q, r)$ where $r|q-1$.
\qed

\bigskip

\bigskip

Now we  prove our main theorem, that is,  determine  the family of $2$-distance-transitive  Cayley graphs over the generalized quaternion  groups.

\bigskip

\noindent {\bf Proof of Theorem \ref{2dt-gq-theo1}.}
Let $\Gamma$ be a connected $(G,2)$-distance-transitive Cayley graph over the group $T=\langle a,b|a^{2n}=1,a^n=b^2,a^b=a^{-1}\rangle$, where $n\geq 2$
and $T\leq G:=\Aut(\Gamma)$.
Then $\Gamma$ has at least 8 vertices.
By Remark \ref{rem-1},    $\Gamma$ has valency at least 4.

Suppose that  $G$ is quasiprimitive on $V(\Gamma)$. Then it follows from
Lemma  \ref{2dt-cay-quasi-th1} that
$\Gamma$ is the complete graph $\K_{4n}$, contradicting that $\Gamma$ has diameter at least 2.

Thus    $G$ is not quasiprimitive on $V(\Gamma)$, and so
 $G$ has at least one intransitive normal subgroup, and this subgroup has at least two orbits on $V(\Gamma)$.

If  every non-trivial normal subgroup of $G$ has at most two orbits on $V(\Gamma)$ and there exists one which has exactly two orbits on $V(\Gamma)$, then $G$ is bi-quasiprimitive on $V(\Gamma)$ and  $\Gamma$ is a bipartite graph. Then by Lemma \ref{2dt-cay-biquasi-th1},
 $\Gamma$ is one of the following graphs:
$\K_{2n,2n}$,  $ \K_{2n,2n}-2n\K_2$ ($n\geq 3$),  $B(\PG(d-1,q))$ and $B'(\PG(d-1,q))$, where $d\geq 3$, $q$ is a prime power and $\frac{q^d-1}{q-1}=2n$.

Now let $N$ be a  normal subgroup of $G$ maximal with respect to having at least $3$ orbits.
Let $H=\langle a\rangle$. Then $H$ has 2 orbits on $V(\Gamma)$, and so $\Gamma$ is a bicirculant.
 Suppose that $\Gamma\ncong \K_{m[f]}$ for any  $m\geq 3$ and $f\geq 2$.
  Then by Theorem \ref{bic-reduction-th2},
$\Gamma$ is a cover of $\Gamma_N$ which
 is $(G/N,s)$-distance-transitive, where $s=\min\{2,\diam(\Gamma_N)\}$, and $G/N$ is faithful and either quasiprimitive or bi-quasiprimitive on $V(\Gamma_N)$.

Moreover,   $N=\langle a^i\rangle$ where $i$ is a proper divisor of $2n$,  and one of the following holds:

\begin{enumerate}[$(1)$]
\item    $T/N\cong \mathbb{Z}_2\times \mathbb{Z}_2$, $i=2$ and $n$ is even;

\item    $T/N\cong \mathbb{Z}_4$, $i=2$ and $n$ is odd;

\item      $T/N\cong Q_{2i}$, $i\geq 3$ is even  and $i$ is not a divisor of $n$;

\item      $T/N\cong D_{2i}$, $i\geq 3$  and $i$ is  a divisor of $n$.
\end{enumerate}

If case (1) or (2) holds, then $\Gamma_N$ has 4 vertices, and so $\Gamma_N$ is $\K_4$ or $C_4$, contradicting that $\Gamma$ has valency at least 4.



For the     case (4),   Lemma \ref{2dt-cay-qd-lem-1} says that
 $\Gamma$ is one of the following graphs:
 \begin{itemize}

\item[(1)]        $\K_{2n,2n}-2n\K_2$ with $n\geq 3$;

\item[(2)]        $ X_1(4,q)$ where $q\equiv 3\pmod{4}$ and  $q=n-1$;

\item[(3)]  $\Gamma(d, q, r)$ where $r|q-1$;

\item[(4)] $X(2,2)$;

\item[(5)] $X'(3,2)$.

\end{itemize}

It remains to consider   case (3). In this case, $T/N\cong Q_{2i}$, where $i\geq 3$ is even  and $i$ is not a divisor of $n$.
If $G/N$ is  bi-quasiprimitive on $V(\Gamma_N)$, then applying Lemma \ref{2dt-cay-qd-lem-2},
$\Gamma$ is one of the following graphs:
 \begin{itemize}


\item[(1)]   $\K_{q+1}^{2d}$, where  $q$ is an odd prime power, $d\geq 2$   dividing $q-1$   and $2n=d(q+1)$;

\item[(2)] $\Gamma(d, q, r)$ where $r|q-1$;



\item[(3)]  $ X_2(3)$;

\item[(4)] $X(2,2)$;

\item[(5)] $X'(3,2)$.

\end{itemize}

From now on, we assume that    $G/N$ is  quasiprimitive  on $V(\Gamma_N)$. Then by Lemma \ref{2dt-cay-quasi-th1},
$G/N$ acts  2-transitively on $V(\Gamma_N)$ and  $\Gamma_N$ is a complete graph, say $\K_r$.
Furthermore, as $N\leq H$, each $N$-orbit is contained in either $H_0$ or $H_1$, and so
$H/N$ has 2 orbits on $V(\Gamma_N)$, as a result $\Gamma_N$ is a bicirculant over the cyclic group $H/N$, and so $r$ is even.

 If $G/N$ acts $3$-transitively on $V(\Gamma_N)$, then
 $\Gamma_N$ is  $(G/N,2)$-arc-transitive. By Lemma \ref{2at-2at-1},
$\Gamma$ is   $(G,2)$-arc-transitive, and further,
 \cite[Lemma 4.1]{WJ-2023jctb} indicates that  $\Gamma$ is either

 \begin{itemize}

\item[(1)] $\K_{2n,2n}-2n\K_2$ with $n\geq 3$; or

\item[(2)]    $ X_1(4,q)$ where $q\equiv 3\pmod{4}$ and  $q=n-1$.
 \end{itemize}

Now we consider the case   that  $\Gamma_N$ is not $(G/N,2)$-arc-transitive. Then  $G/N$ is $2$-transitive but not $3$-transitive on $V(\Gamma_N)$.
Note that $G/N$ has a cyclic subgroup $H/N$ which has 2 orbits on $V(\Gamma_N)$.
Hence $G/N$ is listed as in Theorems \ref{bicirculant-primitive-2} and \ref{2-trans-group}.

By Proposition \ref{3trans-infint-1},
 one of the following cases of  Theorem \ref{2-trans-group} occurs:

 \begin{itemize}

\item[(1)(b)]  $r=(q^d-1)/(q-1)$  and $\PGL(d,q)\leqslant G/N\leqslant  P \Gamma L(d,q)$ for some odd  prime power $q$ and $d\geqslant 3$ even.

\item[(2)(a)]    $G/N\leq A\Gamma L(1,16)$, the group of affine transformations of $F_{16}$ extended by the automorphisms of $F_{16}$.

\item[(2)(b)]  $SL(4,2)\unlhd (G/N)_0$, $r=16$.

\item[(2)(e)]  $(G/N)_0\cong A_6$ or $A_7$, $r=16$.

 \end{itemize}


Applying   Lemma \ref{2dt-volt-lem1}, we can assume  that the covering transformation group
$N$ is a cyclic group $\mathbb{Z}_p$ for some prime $p$.

For cases (2)(a), (2)(b) and  (2)(e), we have $Q: = soc(G/N) \cong  \mathbb{Z}_4^2$
and $\Gamma_N\cong \K_{16}$. Let $Y$ be a normal subgroup of $G$ such that  $N\leq Y$ and
$Y/N = Q$. Then $Y = N.Q$ is transitive on $V (\Gamma)$. Assume first that $p = 2$. Then $\Gamma$ is a $(G,2)$-distance-transitive $\mathbb{Z}_2$-cover of $\Gamma_N\cong  \K_{16}$, and hence $\Gamma$ has order 32 and valency 15.
It follows from \cite{Conder-1} that
there are  3 such candidates, and by Magma \cite{Magma-1997},
two of them are not bicirculants, and
the other one is isomorphic to $\K_{16,16} - 16\K_2$. However, every  arc-transitive subgroup
of $\Aut(\K_{16,16}-16\K_2)$ containing $Y$ as a normal regular subgroup has no regular subgroup
isomorphic to $T\cong Q_{32}$, which is a contradiction.

Assume now that $p$ is an odd prime. Let $Y_2$ be a Sylow 2-subgroup of $Y$. Then $Y_2\cong  \mathbb{Z}_4^2$, and
$Y = N : Y_2\cong  \mathbb{Z}_p.\mathbb{Z}_4^2$. Set  $C: = C_Y (N)$. Then $C = N \times  C_2$, where $C_2$ is a Sylow 2-subgroup of $C$. Applying the $N/C$-theorem, $Y/C$ is isomorphic
to a subgroup of $\Aut(N) \cong  \mathbb{Z}_{p-1}$, which implies that $Y/C$ is cyclic. As $Y_2\cong \mathbb{Z}_2^4$, we
 conclude that 8 is a divisor of $|C|$ and hence $C_2\cong  \mathbb{Z}_3^2$ or $\mathbb{Z}_4^2$. Since $Y $ is a normal subgroup of $ G$, we have $C \unlhd G$, and due to $C_2$ is
characteristic in $C$,  $C_2$ is a normal subgroup of $ G$. Note that   $C_2$ has at least 3 orbits on $V (\Gamma)$. If $\Gamma$ is
a cover of the normal quotient graph $\Gamma_{C_2}$, then by Theorem \ref{bic-reduction-th2}, we would have  $C_2 \leq  H$,
which is   impossible as $H$ is cyclic and $C_2\cong \mathbb{Z}_3^2$ or $\mathbb{Z}_4^2$. Thus, $\Gamma$ is not
 a cover of $\Gamma_{C_2}$. Since $Y$ acts   transitively  on the vertex set of  $\Gamma$, it follows from \cite[Lemma 5.3]{DGLP-locdt-2012}    that
$\Gamma\cong \K_{m[b]}$
for some integers $m\geq  3$ and $b\geq  2$ with $mb = 4n$, a contradiction.

Now we consider the case    (1)(b). Then  $G/N = PSL(2, q).o$ contains the subgroup $T/N $.
Note that   $N=\langle a^i\rangle$ is  a normal subgroup of $T\cong Q_{4n}$,    where $i$ is a divisor of $2n$ and $i\neq 1,2n$, and $|N|=2n/i=p$.

If $i=2$, then $T/N\cong \mathbb{Z}_2\times \mathbb{Z}_2$ or $\mathbb{Z}_4$ according to whether $n$ is even or odd, respectively.
Then $\Gamma_N$ has 4 vertices, as a result   $\Gamma_N\cong \K_4$, and so $\Gamma$ has valency 3.
By Remark \ref{rem-1}, this case does not occur.

Assume  $i\geq 3$. If $i$ is not a divisor of $n$, then $i$ is even, and $T/N\cong Q_{4(i/2)}$; and
if $i$ is  a divisor of $n$, then  $T/N\cong D_{2i}$.

Suppose first that $i$ is  a divisor of $n$ and    $T/N\cong D_{2i}$.
Then by \cite[Theorem 3.3]{SLZ-2014},
$d=2, r = q + 1$ and $G/N = PSL(2, q).o$, where $q = e^f \equiv 3 \pmod  4$ with $e$ a prime,
and $o\leq \mathbb{Z}_2\times \mathbb{Z}_f$ does not contain the diagonal automorphism of $PSL(2, q)$.

Thus  $T/N \leq  PSL(2, q)$, and so $T\cong  D_{p(q+1)}$ as $N\cong \mathbb{Z}_p$. Let $Q :=
soc(G/N)$. Since $G/N = Q.o$, there is a normal subgroup $Y$ of $ G$ such that $Y/N = Q$. As
$N\cong \mathbb{Z}_p$ with $p$ a prime, it follows from  \cite[Lemma 2.11]{PLHL-2014}  that $Y = N.Q$ is a central extension,
$Y = NY'$, and $Y'\cong  M.Q$ with $M\leq  Mult(Q) \cap N$.

Since $q \equiv 3 \pmod 4$, it follows from \cite[P. 302, Table 4.1]{Gorenstein-1} that $Mult(PSL(2, q)) \cong \mathbb{Z}_2$.
If $p  = 2$, then \cite[Section 3.3.6]{Wilson-book} implies that $Y = \mathbb{Z}_2\times   PSL(2, q)$ or $SL(2, q)$.
Note that $D_{2(q+1)}\cong   T \leq   Y$ as $T/N \leq Q = Y/N$. Then $Y$ has an element of
order $q + 1$, say $z$. By \cite[Lemma 2.9]{Pan-2014}, $SL(2, q)$ has no subgroup isomorphic to $D_{2(q+1)}$,
and we may let $Y = \mathbb{Z}_2 \times  PSL(2, q)$. Then $z = z_1z_2$ with $z_1\in  \mathbb{Z}_2$ and $z_2\in  PSL(2, q)$,
and since $q + 1$ is even, $z_2$ has order $q + 1$, which is impossible by \cite[Table 8.1]{BHR-2013}. Thus,
$p$ is an odd prime, and so $Mult(Q) \cap  N = 1$, forcing   $M = 1$. Then $Y = N \times   Y'$,
and hence $N\leq  Z(Y )$, the center of $Y$. This is impossible as   $N < T \leq  Y$ and
$T\cong  D_{p(q+1)}$.

In the remainder we suppose that   $i> 3$  is not a divisor of $n$,  $i$ is even, and $T/N\cong Q_{4(i/2)}$.

Let $Y$ be a minimal subgroup of $G/N$ which contains $T/N$ and is primitive on $V(\Gamma_N)$.
Then  by Lemma \ref{2dt-cay-quasi-th1},  $Y$ is 2-transitive, and so $soc(Y )$ is primitive. Let $K = Y_u$, where $u\in V(\Gamma_N)$, and let $M$ be a maximal subgroup of $Y$
which contains $T/N$. Then, $Y = MK$. By the assumption on the minimality of $Y$,
$M$ is imprimitive on $V(\Gamma_N)$, and we conclude that $soc(Y ) \leq  M$. Hence, $Y = MK$ is a
maximal factorization of $Y$, which is given in \cite{LPS-1990}. Now $Y$ is 2-transitive on $[Y : K]$
of degree $r$, and $r = |T/N|$ is even, say $r = 2l$.

Suppose that $soc(Y ) = PSL(d, q)$ and $2l = \frac{q^d-1}{q-1}$ with $d\geq  3$. Then, the maximal
factorization $Y = MK$ is given in \cite[Table 1]{LPS-1990}, such that $K = P_1$ is a parabolic
subgroup. By \cite[Table 1]{LPS-1990}, either $d$ is even and $soc(M) \cong  PSp(d, q)$, or $M = \tensor*[^\wedge]{G}{}  L(a, q^b).\mathbb{Z}_b$ where $ab = d$. If $soc(M) = PSp(d, q)$, then it follows that $M$ is
quasiprimitive, and so $M$ is 2-transitive of degree $\frac{q^d-1}{q-1}$, which is not
possible. Thus, we have that $M = \tensor*[^\wedge]{G}{} L(a, q^b).\mathbb{Z}_b$. Suppose that $a = 1$. Then
$M/Z(M) = PGL(a, q^b).\mathbb{Z}_b$. Since   $T/NZ(M)/Z(M)\cong
(T/N)/(T/N\cap  Z(M))$ has order $\frac{q^{ab}-1}{q^b-1}$, we have   $|T/N \cap  Z(M)| = \frac{q^b-1}{q-1}$. Hence  $T/N$ is a central extension of $T/N\cap Z(M)$ by $(T/N)/(T/N\cap  Z(M))$, which is a contradiction as  $\frac{q^b-1}{q-1}=|G\cap  Z(M)|\leq |Z(G)| \leq  2$. Consequently   $a = 1$, and so $M = P\Gamma L(1, q^d)$. However, it is easy to  show that $P\Gamma L(1, q^d)$ with $d\geq  3$ contains no generalized quaternion  subgroup of order
$\frac{2(q^d-1)}{q-1}$, which is again a contradiction

This concludes   the proof.
 \qed

\end{document}